\documentclass[11pt]{article}
\usepackage{latexsym,amsmath,color,amsthm,amssymb,epsfig,graphicx,mathrsfs}
\usepackage{graphicx}
\usepackage{fullpage}
\def\qed{\hfill \ifhmode\unskip\nobreak\fi\quad\ifmmode\Box\else$\Box$\fi\\ }

\newtheorem{thrm}{Theorem}[section]
\newtheorem{lemm}[thrm]{Lemma}
\newtheorem{propo}[thrm]{Proposition}
\newtheorem{coro}[thrm]{Corollary}
\newtheorem{defi}[thrm]{Definition}
\newtheorem{conjecture}[thrm]{Conjecture}

\newtheorem{con}{Construction}

\newcommand{\thm}{\begin{thrm}}
\newcommand{\xthm}{\end{thrm}}
\newcommand{\lem}{\begin{lemm}}
\newcommand{\xlem}{\end{lemm}}
\newcommand{\prf}{\begin{proof}}
\newcommand{\xprf}{\end{proof}}
\newcommand{\prop}{\begin{propo}}
\newcommand{\xprop}{\end{propo}}
\newcommand{\cor}{\begin{coro}}
\newcommand{\xcor}{\end{coro}}
\newcommand{\defn}{\begin{defi}}
\newcommand{\xdefn}{\end{defi}}
\newcommand{\conj}{\begin{conjecture}}
\newcommand{\xconj}{\end{conjecture}}

\renewcommand{\phi}{\varphi}

\newtheorem{theorem}{Theorem}

\newtheorem{lemma}[theorem]{Lemma}
\newtheorem{proposition}[theorem]{Proposition}

\begin{document}

\title{\vspace{-0.5in} The structure  of large intersecting families}

\author{
{\large{Alexandr Kostochka}}\thanks{
\footnotesize {University of Illinois at Urbana--Champaign, Urbana, IL 61801
 and Sobolev Institute of Mathematics, Novosibirsk 630090, Russia. E-mail: \texttt {kostochk@math.uiuc.edu}.
 Research of this author
is supported in part by NSF grant  DMS-1266016
and and by
grants 15-01-05867 and 16-01-00499  of the Russian Foundation for Basic Research.
}}
\and
{\large{Dhruv Mubayi}}\thanks{
\footnotesize {Department of Mathematics, Statistics, and Computer Science, University of Illinois at Chicago, Chicago, IL 60607.
E-mail:  \texttt{mubayi@uic.edu.}
Research partially supported by NSF grant DMS-1300138.}}}
\maketitle

\vspace{-0.3in}

\begin{abstract}
A collection of sets is {\em intersecting} if every two members have nonempty intersection. We describe the structure of  intersecting families of $r$-sets of an $n$-set whose size is quite a bit smaller than the maximum  ${n-1 \choose r-1}$ 
given by the Erd\H os-Ko-Rado Theorem.  In particular, this extends the Hilton-Milner theorem on nontrivial intersecting families and answers a recent question of Han and Kohayakawa for large $n$. 
In the case $r=3$ we describe the structure of all intersecting families with more than 10 edges.
We also prove a stability result for the Erd\H os matching problem.  Our short proofs are simple applications of 
the Delta-system method introduced and  extensively used by Frankl since 1977.
\end{abstract}

\section{Introduction}
An $r$-uniform hypergraph $H$, or simply {\em $r$-graph}, is a family of $r$-element subsets of a finite set.
We associate an $r$-graph $H$ with its edge set and call its vertex set $V(H)$. Say that $H$ is {\em intersecting} if $A \cap B \ne \emptyset$ for all $A, B \in F$. 
A {\em matching} in $H$ is a collection of pairwise disjoint sets from $H$. A {\em vertex cover} (henceforth cover)  of $H$ is a set of vertices 
intersecting every edge of $H$.  Write $\nu(H)$ for the size of a maximum matching and
 $\tau(H)$ for the size of a minimum  cover of $H$. Say that $H$ is {\em trivial} or a {\em star} if $\tau(H)=1$, otherwise call $H$ {\em nontrivial}.
 
A fundamental problem in the extremal theory of finite sets is to determine the maximum size of an $n$-vertex $r$-graph $H$  with $\nu(H) \le s$. 
The case $s=1$ is when $H$ is intersecting, and in this case the Erd\H os-Ko-Rado Theorem~\cite{EKR} states that the maximum is 
${n-1 \choose r-1}$ for $n \ge 2r$ and if $n>2r$, then equality holds only if $\tau(H)=1$. 
More generally,   Erd\H os~\cite{E1965} proved the following.

\begin{theorem} [Erd\H os~\cite{E1965}] \label{em} For $r \ge 2$, $s\geq 1$ and $n$ sufficiently large, every
 $n$-vertex $r$-graph $H$  with $\nu(H) \le s$, satisfies
\begin{equation}\label{1125}
|H|\leq em(n,r,s):={n \choose r} - {n-s \choose r} \sim s{n \choose r-1},
\end{equation}
and if equality in~\eqref{1125} holds, then $H$ is the $r$-graph $EM(n,r,s)$ described below.
\end{theorem}

\begin{con}\label{con2}
Let $EM(n,r,s)$ be the $n$-vertex $r$-graph that has $s$ special vertices $x_1,\ldots,x_s$ and the edge set
consists of the all $r$-sets intersecting $\{x_1,\ldots,x_s\}$. In particular, $EM(n,r,1)$ is a full star.
\end{con}

There has been a lot of recent activity  on Theorem~\ref{em} for small $n$ (see, e.g., \cite{Frankl2013, Frankl3uniform, HLS, LM}). 

 Hilton and Milner~\cite{HM} proved a strong stability result for the Erd\H os-Ko-Rado Theorem:

\begin{theorem}[Hilton-Milner \cite{HM}, Proposition $\mathcal{T}$] \label{hmt}
Suppose that $ 2\leq r\leq n/2$ and $|H|$ is an $n$-vertex intersecting $r$-graph with $\tau(H)\geq 2$.
Then 
\begin{equation}\label{ehm}
 |H|\leq hm(n,r) := {n-1 \choose r-1}-{n-r-1 \choose r-1}+1 \sim r{n \choose r-2},
\end{equation}
and if equality in~\eqref{ehm} holds, then $H$ is the $r$-graph $HM(n,r)$ described below.
\end{theorem}

\begin{con}\label{con1}
  For $n \ge 2r$, let
 $HM(n,r)$ be the following $r$-graph on $n$ vertices: Choose an $r$-set $X=\{x_1, \ldots, x_r\}$ and 
a   special vertex $x \not\in X$, and let $HM(n,r)$ consist of the set $X$ and all $r$-sets containing $x$ and a vertex of $X$. 
\end{con}

Observe that $HM(n,r)$ is intersecting, $\tau(HM(n,r))=2$, and $|HM(n,r)|=hm(n,r)$.
Bollob\' as,  Daykin and Erd\H os~\cite{BDE} extended Theorem~\ref{hmt} to $r$-graphs with matching number $s$ in the way
Theorem~\ref{em} extends the Erd\H os-Ko-Rado Theorem.

\begin{theorem}[Bollob\' as-Daykin-Erd\H os~\cite{BDE}, Theorem 1] \label{bde}
Suppose  $ r\geq 2$, $s\geq 1$ and $n> 2r^3s$. If
 $H$ is an $n$-vertex $r$-graph with $\nu(H)\leq s$ and $|H|>em(n,r,s-1)+  hm(n,r)$,
then $H\subseteq EM(n,r,s)$. 
\end{theorem}

Han and  Kohayakawa~\cite{HK} refined  Theorem~\ref{hmt} using the following construction.

\begin{con}\label{con3}
For $r \ge 3$, the $n$-vertex $r$-graph $HM'(n,r)$ has $r+2$ distinct special vertices $x, x_1, \ldots, x_{r-1}, y_1, y_2$ and all edges $e$ such that\\
1) $\{x, x_i\} \subset e$ for any $i \in [r-1]$, or\\
2) $\{x, y_1, y_2\} \subset e$, or\\
3) $e=\{x_1, \ldots, x_{r-1}, y_1\}$, or $e=\{x_1, \ldots, x_{r-1}, y_2\}$.
\end{con}

Note that $HM'(n,r)$ is intersecting, $\tau(HM'(n,r)) =2$, and
$HM'(n,r) \not\subset HM(n,r)$. Let $hm'(n,r)=|HM'(n,r)|$ so that
$$hm'(n,r)= {n-1 \choose r-1} - {n-r \choose r-1} +{n-r-2 \choose r-3} +2 \sim (r-1){n \choose r-2}.$$ 
The  result of~\cite{HK} for $r\geq 5$ is:

\begin{theorem}[Han-Kohayakawa~\cite{HK}]\label{th2} Ler $r \ge 5$ and $n>2r$. 
If $H$ is an $n$-vertex intersecting $r$-graph, $\tau(H) \ge 2$ and $|H|\geq hm'(n,r)$, then $H \subseteq HM(n,r)$ or
 $H=HM'(n,r)$.
\end{theorem}
They also resolved the cases $r=4$ and $r=3$, where the statements are similar but somewhat more involved.

For large $n$ Frankl~\cite{Frankl1980} gave an exact upper bound on the size of intersecting $n$-vertex $r$-graphs $H$
with $\tau(H)\geq 3$. He introduced the following family.  We write $A+a$ to mean $A \cup \{a\}$.

\begin{con}[\cite{Frankl1980}]\label{con8} The vertex set $[n]$ of the $n$-vertex $r$-graph $FP(n,r)$ 
contains a special subset $X=\{x\}\cup Y\cup Z$ with $|X|=2r$ such that
$|Y|=r$, $|Z|=r-1$, where a subset $Y_0=\{y_1,y_2\}$ of $Y$ is specified.
The edge set of $FP(n,r)$ consists of all $r$-subsets of $[n]$ containing a member of the 
family
$$G=\{A\subset X\,: \,|A|=3, x\in A, A\cap Y\neq\emptyset, A\cap Z\neq\emptyset\}
\cup \{Y, Y_0+x, Z+y_1,Z+y_2\}.
$$
\end{con}

By construction, $FP(n,r)$ is an intersecting $r$-graph with $\tau(FP(n,r))=3$. Frankl proved the following.

\begin{theorem}[Frankl~\cite{Frankl1980}] \label{fp} Let  $r\geq 3$ and $n$ be sufficiently large.
Then every   intersecting $n$-vertex
$r$-graph $H$ with $\tau(H)\geq 3$ satisfies $|H|\leq |FP(n,r)|$. Moreover, if $r\geq 4$, then
equality is attained only if $H=FP(n,r)$.
\end{theorem}
He used the following folklore result.

\begin{proposition} \label{folk} Every intersecting $3$-graph $H$ with $\tau(H)\geq 3$ satisfies $|H|\leq 10$.
\end{proposition}
Note that Erd\H os and Lov\' asz~\cite{EL} proved the more general result that for every $r\geq 2$ each intersecting 
$r$-graph $H$ with $\tau(H)=r$ has at most $r^r$ edges. But their proof gives the  bound $25$ for $r=3$, while
 Proposition~\ref{folk} gives $10$.

 In this short paper, we determine for large $n$, the structure of $H$ in the situations described  above when $|H|$ is somewhat smaller than  the bounds in 
 Theorems~\ref{th2} and~\ref{ehm}. In particular, our Theorem~\ref{th3} below answers
for large $n$ the question of Han and  Kohayakawa~\cite{HK} at the end of their paper.
We also use  Theorem~\ref{fp} to describe large dense hypergraphs $H$ with $\nu(H)\leq s$ and $\tau(H)=2$.
Related results  can be found in~\cite{Frankl1980, Frankl1987}.

 \section{Results}


First we characterize the 
nontrivial intersecting $r$-graphs that have a bit fewer edges than $hm'(n,r)$.  We need to describe three constructions before we can state our  result.

\begin{con}\label{con4}
For $r \ge 3$, $1\leq t\leq r-1$ and $t=n-r$, the  $n$-vertex $r$-graph $HM(n,r,t)$ has $r+t$
 distinct special vertices $x, x_1, \ldots, x_{r-1}, y_1, y_2, \ldots, y_t$ and all edges $e$ such that\\
1) $\{x, x_i\} \subset e$ for any $i \in [r-1]$, or\\
2) $e=\{x_1, \ldots, x_{r-1}, y_j\}$ for all $1\leq j\leq t$, and\\
3) if $1\leq t \le r-1$, then  include all $e$ such that $\{x, y_1, \ldots, y_t\} \subseteq e$.  
\end{con}

 Let $hm(n,r,t)=|HM(n,r,t)|$. 
 Note that $HM(n,r,1)=HM(n,r)$, and $HM(n,r,2)=HM'(n,r)$. For $n$ large, we have the inequalities
 $$hm(n,r)=hm(n,r,1)>\cdots > hm(n,r,r-1)=hm(n,r,r)  < hm(n,r,n-r).$$
Note that $HM(n,r,t)$ is intersecting, $\tau(HM(n,r,t)) =2$, and
$HM(n,r,t) \not\subseteq HM(n,r,t-1)$. Also, for fixed $r \ge 4$ and $2 \le t \le n-r$,
$$hm(n,r,t) \sim (r-1) {n \choose r-2}.$$

\begin{con}\label{con5}
The $n$-vertex $r$-graph $HM(n,r,0)$ has $3$ special vertices $x, x_1, x_2$ and all
edges that contain at least two of these $3$ vertices.
\end{con}

 By definition,
\begin{equation}\label{112}
 |HM(n,r,0)|=3{n-3\choose r-2}+{n-3\choose r-3}. 
\end{equation}



\begin{con}\label{con6}
The  $n$-vertex $r$-graph $HM''(n,r)$ has $r+4$ special vertices $x, x_1, \ldots, x_{r-2}$ and $y_1, y_1', y_2, y_2'$
 and all edges $e$ such that\\
1) $\{x, x_i\} \subset e$ for some $i \in [r-2]$, or\\
2) $\{x, y_1, y_2\} \subseteq e$, or $\{x, y_1, y_2'\} \subseteq e$ or  $\{x, y_1', y_2\} \subseteq e$ or $\{x, y_1', y_2'\} \subseteq e$, or\\
3) $e=\{x_1, \ldots, x_{r-2}, y_1, y_1'\}$, or
 $e=\{x_1, \ldots, x_{r-2}, y_2, y_2'\}$.
\end{con}

Note that $HM''(n,r)$ is intersecting, $\tau(HM''(n,r)) =2$, and
$HM''(n,r) \not\subseteq HM(n,r,t)$ for any $t$. Let $hm''(n,r)=|HM''(n,r)|$ so that for $ r \ge 5$, 
\begin{align}\label{1122}
hm''(n,r)&={n-1 \choose r-1} - {n-r+1 \choose r-1} +4{n-r-3 \choose r-3} +4{n-r-3 \choose r-4}+{n-r-3 \choose r-5}+
2 \notag \\
&\sim (r-2){n \choose r-2}.
\end{align}

\begin{theorem}\label{th3} Fix $r \ge 4$. Let $n$ be   sufficiently large.
 If $H$ is an $n$-vertex intersecting $r$-graph with $\tau(H) \ge 2$ and $|H|>hm''(n,r)$, 
then $H \subseteq HM(n,r,t)$ for some $t\in\{1,\ldots,r-1,n-r\}$ or $r=4$ and $H \subseteq HM(n,4,0)$.
The bound on $H$ is sharp due to $HM''(n,r)$.
\end{theorem}

When $r=3$ we are able to obtain stronger results than Theorem~\ref{th3}, and  describe the structure of {\em 
almost all} intersecting $3$-graphs. We will use the following construction.

\begin{con}\label{con7} Let $n\geq 6$.\\
$\bullet$ For $i=0,1,2$, let 
$$H_i(n)=HM(n,3,i) \qquad \hbox{\rm and} \qquad   H(n)=EM(n,3,1).$$

$\bullet$   The $n$-vertex $3$-graph $H_3(n)$  has special vertices $v_1,v_2,y_1,y_2,y_3$ and 
its edges are the $n-2$ edges containing $\{v_1,v_2\}$ and the $6$ edges each of which contains one
of $v_1,v_2$ and two of
$y_1,y_2, y_3$.

$\bullet$ Each of the $n$-vertex $3$-graphs $H_4(n)$ and $H_5(n)$
 has $6$ special vertices $v_1,v_2,z_{1,1}z_{1,1}',z_{2,1}z_{2,1}'$ and contains
 all edges containing $\{v_1,v_2\}$. Apart from these, $H_4(n)$ contains edges
 $$v_1z_{1,1}z_{1,1}', v_1z_{2,1}z_{2,1}', v_2z_{1,1}z_{2,1},  v_2z_{1,1}z'_{2,1},  v_2z'_{1,1}z_{2,1},  v_2z'_{1,1}z'_{2,1}$$
and $H_5(n)$ contains edges
 $$v_1z_{1,1}z_{1,1}', v_1z_{2,1}z_{2,1}', v_1z_{1,1}z'_{2,1},  v_2z_{1,1}z'_{2,1},  v_2z_{1,1}z_{2,1},  v_2z'_{1,1}z'_{2,1}.$$
\end{con}


\begin{theorem} \label{th4} Let  $H$ be an  intersecting
$3$-graph and $n=|V(H)|\geq 6$. If $\tau(H)\leq 2$, then
 $H$ is contained in one of $H(n), H_0(n),  \ldots,H_5(n)$.
This yields that\\
(a) if  $|H|\geq 11$, then  $H$ is contained in one of $H(n), H_0(n),  \ldots,H_5(n)$;\\
(b) if  $|H|>n+4$, then $H$ is contained in $H(n),H_0(n), H_1(n)$ or $H_2(n)$.
\end{theorem}

The restriction $|H|\geq 11$ cannot be weakened because of $K_5^3$ and $|H|>n+4$ cannot be weakened because $|H_3(n)|= |H_4(n)|=|H_5(n)|=n+4$.

To prove an analog of Theorem~\ref{th4} for $r$-graphs, we need an  extension of Construction~\ref{con7}:

\begin{con}\label{con9} Let $n\geq r+1$. For $i=0,\ldots,5$, let the  $r$-graph
$H^r_i(n)$ have the vertex set of the $3$-graph $H_i(n)$ and the edge set of $H^r_i(n)$
consist of all $r$-tuples containing an edge of $H_i(n)$.
\end{con}
By definition, $H_0^r(n)=HM(n,r,0)$. Each $H^r_i(n)$ is intersecting, since each $H_i(n)$
is intersecting. Using Theorem~\ref{fp}, we extend Theorem~\ref{th4} as follows:

\begin{theorem} \label{th4'} Let  $r\geq 4$ be fixed and $n$ be sufficiently large.
Then there is $C>0$ such that for every
  intersecting $n$-vertex
$r$-graph $H$ with $|H|> |FP(n,r)|=O(n^{r-3})$, one can delete from $H$ at most $Cn^{r-4}$ edges
so that the resulting $r$-graph $H'$
  is contained in one of $ H^r_0(n),  \ldots,H^r_5(n), EM(n,r,1)$.
\end{theorem}

The results above naturally extend to $r$-graphs $H$ with $\nu(H)\leq s$. For example,
 Theorem~\ref{th3} extends to the  following result which  implies Theorem~\ref{bde} for large $n$.

\begin{theorem} \label{th1}
 Fix $r \ge 4$ and $s\geq 1$. Let $n$ be   sufficiently large.
 If $H$ is an $n$-vertex  $r$-graph with $\nu(H)\leq s$ and $|H|>em(n,r,s-1)+hm''(n-s+1,r)$, 
then $V(H)$ contains a subset $Z=\{z_1,\ldots,z_{s-1}\}$ such that either $\tau(H-Z)=1$ or
$H-Z \subseteq HM(n-s+1,r,t)$ for some $t\in\{1,\ldots,r-1,n-s+1-r\}$ or $r=4$ and $H-Z \subseteq HM(n-s+1,4,0)$.
 The bound on $|H|$ is sharp.
\end{theorem}

Theorems~\ref{th2} and~\ref{th4'} can be extended in a similar way. We leave this to the reader.

\section{Proof of Theorem~\ref{th3}}
The main tool used in the proof is the Delta-system method developed by Frankl (see, e.g.~\cite{Frankl1977,Frankl1980}). 
Recall that a  $k$-{\em sunflower} $S$ is a collection of distinct sets $S_1, \ldots, S_k$ such that for every $1\le i<j \le k$, we have $S_i \cap S_j =\bigcap_{\ell=1}^k S_{\ell}$. The common intersection of the $S_i$ is the {\em core} of $S$.  We will use the following fundamental result of Erd\H os and Rado~\cite{ER}.

\begin{lemma} [Erd\H os-Rado Sunflower Lemma \cite{ER}]
\label{er} 
For every $k,r\ge 2$ there exists  $f(k,r)< k^rr!$ such that the following holds: every $r$-graph $H$ with no $k$-sunflower satisfies $|H|<f(k,r)$.
\end{lemma}

{\bf Proof of Theorem \ref{th3}.} Let $r\geq 4$ and  $H$ be an $n$-vertex intersecting $r$-graph with $\tau(H) \ge 2$ and $|H|>hm''(n,r)$.
 Define $B^*(H)$ to be the set of $T \subset V(H)$ such that 

(i)  $0<|T|<r$, and 

(ii) $T$ is the core of an  $(r+1)^{|T|}$-sunflower in $H$.

Define $$B'(H)=\{T \in B^*(H): \nexists U \in B^*(H), U \subsetneq T\}$$ to be the set of all inclusion minimal elements in $B^*(H)$. Next, 
let 
$$B''(H)=\{e \in H: \nexists T \subsetneq e, T \in B^*(H)\}$$ be the set of edges in $H$ that contain no member of $B^*(H)$. 
Finally, set 
$$B(H)=B'(H) \cup B''(H).$$
 Let $B_i$ be the sets in $B(H)$ of size $i$.  
Note that $B_1=\emptyset$ for otherwise we have an $(r+1)$-sunflower with core of size one and since $H$ is intersecting, 
this forces $H$ to be trivial. The following crucial claim proved by Frankl can be found in Lemma 1 in~\cite{Frankl1977, Frankl1980}.
 
 {\bf Claim.} $B_i$ contains no $(r+1)^{i-1}$-sunflower.
 
{\bf Proof of Claim.} Suppose for contradiction that
$S_1, \ldots, S_{(r+1)^{i-1}}$ is an $(r+1)^{i-1}$-sunflower in $B_i$ with core $K$.  
 By definition of $B_i$,  there is an $(r+1)^i$-sunflower ${\cal S}_1=S_{1, 1}, \ldots, S_{1, (r+1)^i}$ in $H$ with core $S_1$. Since $|S_2 \cup \cdots \cup S_{(r+1)^{i-1}}|< (r+1) (r+1)^{i-1}=(r+1)^i$, and ${\cal S}_1$ is an $(r+1)^i$-sunflower, there is a $k=k(1)$ such that 
 $$(S_{1,k(1)} -S_1) \cap (S_2 \cup S_3 \cup \cdots \cup S_{(r+1)^{i-1}}) = \emptyset.$$
Next, we use the same argument to define $S_{2, k(2)}$  such that
$S_{2, k(2)}-S_2$ is disjoint from $S_{1, k(1)} \cup S_3 \cup \cdots\cup  S_{(r+1)^{i-1}}$ and then 
$S_{3, k(3)}$ such that $S_{3, k(3)}-S_3$ is disjoint from $S_{1, k(1)}  \cup S_{2,k(2)} \cup S_3 \cup \cdots \cup S_{(r+1)^{i-1}}$ and so on.
Continuing in this way we finally obtain edges $S_{j, k(j)}$ of $H$
for all $1 \le j \le (r+1)^{i-1}$ that form an $(r+1)^{i-1}$-sunflower  with core $K$.
  This implies that $K \ne \emptyset$ as $H$ is intersecting. Since $|K|\le i-1$, there exists a nonempty $K' \subseteq K$
 such that $K' \in B(H)$. But $K' \subsetneq S_j$ for all $j$, so this contradicts the fact that $S_j \in B(H)$. \qed

Applying the Claim and Lemma~\ref{er} yields $|B_i|< f((r+1)^{i-1}, i)$ for all $i>1$. 
Every edge of $H$ contains an element of $B(H)$ so we can count edges of $H$ by the sets in $B(H)$.  So for $q=|B_2|$  we have 
$$hm''(n,r) < |H| \le \sum_{B \in B_2} {n-2 \choose r-2} + \sum_{i=3}^r \sum_{B \in B_i} {n-i \choose r-i}<q {n-2 \choose r-2} + (r-2)f((r+1)^{r-1},r){n\choose r-3}.$$
Since $hm''(n,r)\sim (r-2){n \choose r-2}$, this  gives $q \ge r-2$. On the other hand, $B_2$ is intersecting and thus the pairs in $B_2$ form either the star $K_{1,q}$ or
a $K_3$.

{\bf Case 1:}  $B_2$ is a $K_3$. Then to keep $H$ intersecting, $H\subseteq HM(n,r,0)$. If $r\geq 5$, then by~\eqref{112} and~\eqref{1122},
$ |HM(n,r,0)|<hm''(n,r)<|H|$, a contradiction. Thus $r=4$ and $H \subseteq HM(n,4,0)$, as claimed.

Since Case 1 is proved, we may assume that $B_2$ is a star  with center $x$ and the set of leaves $X=\{x_1, \ldots, x_{q}\}$.

{\bf Case 2:}  $q\geq r-1$.
 If $ q\ge r$, then $q=r$ and since $H$ is nontrivial,  $H \subseteq HM(n,r)$ and we are done. 
We may therefore assume that $q=r-1$.
Since $\tau(H) \ge 2$, there exists $e$ such that $x \not\in e \in H$, and since $H$ is intersecting we may assume that $e=e_1=X\cup\{ y_1\}$. 
We may also assume that all edges of $H$ that omit $x$ are of the form $e_i=X\cup\{ y_i\}$, where $1 \le i \le t$. 
If $t=1$ then $H \subseteq HM(n,r)$ and we are done, so assume that $t \ge 2$. 
Any edge of $H$ containing $x$ that omits $X$ must contain all $\{y_1, \ldots, y_t\}$. Consequently,
$H \subseteq HM(n,r,t)$ for some $t\in\{1,\ldots,r-1,n-r\}$.

{\bf Case 3:} $q=r-2$.
Let $F_0$ be the set of edges in $H$ that contain $x$ and intersect $X$, $F_1$ be the set of edges of $H$ disjoint from $X$ and
$F_2$ be the set of edges disjoint from $x$. Then $H=F_0\cup F_1\cup F_2$, all edges in $F_1$ contain $x$ and all edges in
$F_2$ contain $X$. Since $|F_0|\le{n-1 \choose r-1} - {n-r+1 \choose r-1}$, by~\eqref{1122},
\begin{equation}\label{1123}
|F_1\cup F_2|>4{n-r-3 \choose r-3} +4{n-r-3 \choose r-4}+{n-r-3 \choose r-5}+
2> 4{n-r-2\choose r-3}.
\end{equation}

Let $G$ be the graph of pairs $ab$ such that $x\notin \{a,b\}$ and $X\cup\{a, b\}\in F_2$. Then $|G|=|F_2|$ and $V(G)\subseteq V(H)-X-\{x\}$.
 
{\bf Case 3.1:}  $\tau(G)=1$. Then $G=K_{1,s}$ for some $1\leq s\leq n-r$. Let the partite sets of $G$ be $x_{r-1}$ and $Y$.
Then every edge in $F_1$ must contain either  $x_{r-1}$ or $Y$. Thus $H \subseteq HM(n,r,t)$ for some $t\in\{1,\ldots,r-1,n-r\}$, as claimed.

{\bf Case 3.2:}  $\tau(G)\geq 2$ and $\nu(G)=1$. Then $G=K_3$ and every edge in $F_1$ must contain at least two vertices of $G$.
Then $|F_1|<3{n-r-1\choose r-3}\sim 3{n\choose r-3}$ and thus $|F_1\cup F_2|=|F_1|+3\sim 3{n\choose r-3}$, contradicting~\eqref{1123}. 

{\bf Case 3.3:}  $\nu(G)\geq 3$. Let $f_1,f_2,f_3$ be disjoint edges in $G$. Then each edge in
$F_1$ has at least $4$ vertices in $f_1\cup f_2\cup f_3\cup \{x\}$ and thus $|F_1|=O(n^{r-4})$. If $F_1=\emptyset$, then
$H\subseteq HM(n,r,n-r)$, as claimed. Suppose there is $e_0\in F_1$. Then each $f\in G$ meets $e_0-x$ and thus
$|G|=|F_2|\leq (r-1)(n-2r+2)+{r-1\choose 2}$. Thus if $r\geq 5$, then $|F_1\cup F_2|\leq O(n^{r-4})+O(n)=o(n^{r-3})$, 
contradicting~\eqref{1123}. Moreover, if $r=4$, then $|F_2|\leq 3(n-6)+3$ and
$|F_1\cup F_2|\leq O(n^{r-4})+3n<4{n-6\choose 1}$, again contradicting~\eqref{1123}.

 {\bf Case 3.4:}  $\nu(G)=2$. Say that a vertex $v$ is {\em big} if $d_G(v)\geq 2r$. Let $v_1,\ldots,v_s$ be all the big vertices in $G$.
 Since $\nu(G)=2$, $s\leq 2$. Since $H$ is intersecting,
 \begin{equation}\label{1124}
\mbox{  Every edge in $F_1$  contains all big vertices.}
\end{equation}
 
Suppose first,  $s=2$. Then to have $\nu(G)=2$, all edges in $F_2$ are incident with $v_1$ or $v_2$; thus $|F_2|<2n$. On the other hand, in this case
 by~\eqref{1124}, $|F_1|\leq {n-r-1\choose r-3}$. Together, this contradicts~\eqref{1123}.
 
 Suppose now, $s=1$. Then to have $\nu(G)=2$, we need $|F_2|\leq d_G(v_1)+2r\leq n+2r$. On the other hand, since $\nu(G)=2$, $G$ has an edge
 $v'v''$ disjoint from $v$. It follows that each edge in $F_1$ meets $v_1v_2$. By this and~\eqref{1123}, $|F_1|\leq 2{n-r\choose r-3}$ and thus
 $|F_1\cup F_2|\leq n+2r+2{n-r\choose r-3}$, contradicting~\eqref{1123}.

Finally, suppose $s=0$. 
Let  edges $y_1y_1'$ and $y_2y_2'$ form a matching in $G$. If $G$ has no other edges, then $H$ is contained in $HM''(n,r)$.
So there is a third edge in $G$. Still, since  $\nu(G)=2$, each edge of $G$ is incident with $\{y_1,y_1',y_2,y_2'\}$ which by
$s=0$ yields $|F_2|=|G|<8r$.
If an edge in $G$ is $y_1y_3$, then each each edge in $F_1$ contains  $\{y_1,y_2\}$ or $\{y_1,y_2'\}$ or
$\{y'_1,y_2,y_3\}$ or $\{y'_1,y'_2,y_3\}$; thus $|F_1|\leq 2{n-r\choose r-3}+2{n-r\choose r-4}\sim 2{n-r\choose r-3}$. This together 
with $|F_2|\leq 8r$ contradicts~\eqref{1123}. If this third edge is $y_1y_2$, then we get a similar contradiction. 
 \qed

\section{On $3$-graphs}

\begin{lemma}\label{le1}
Let $n\geq 6$ and $H$ be an intersecting $3$-graph. If $H$ has a vertex $x$ such that $H-x$ has at most two edges, then
$H$ is contained in one of $H(n), H_0(n),H_1(n), H_2(n), H_4(n)$.
\end{lemma}

{\bf Proof.} If $H-x$ has no edges, then $H\subseteq H(n)$, and if $H-x$ has one edge, then $H\subseteq H_1(n)$.
Suppose $H-x$ has two edges, $e_1$ and $e_2$. If $|e_1\cap e_2|=2$, then we
may assume $e_1=\{x_1,x_2,y_1\}$ and $e_2=\{x_1,x_2,y_2\}$. In this case, each edge in $H-e_1-e_2$ contains $x$ and either
intersects $\{x_1,x_2\}$ or coincides with $\{x,y_1,y_2\}$. This means $H\subseteq H_2(n)$.

 If $|e_1\cap e_2|=1$, then we
may assume $e_1=\{y,v_1,w_1\}$ and $e_2=\{y,v_2,w_2\}$. In this case, each edge in $H-e_1-e_2$ contains $x$ and either
contains $y$ or
intersects each of $\{v_1,w_1\}$ and $\{v_2,w_2\}$. This means $H\subseteq H_4(n)$.\qed

{\bf Proof of of Theorem \ref{th4}.} 
Let  $n\geq 6$ and $H$ be an $n$-vertex intersecting $3$-graph
 with $\tau(H)\leq 2$ not contained in any of $H(n), H_0(n),  \ldots,H_5(n)$.
 Write $H_i$ for $H_i(n)$. 
If $\tau(H)=1$, then $H\subseteq H(n)$.
 So, suppose a set $\{v_1,v_2\}$ covers all edges of $H$, but $H$ is not a star.
Let $E_0=\{e\in H\,:\,
 \{v_1,v_2\}\subset e\}$, and for $i=1,2$, let $E_i=\{e\in H\,:\,
v_{3-i}\notin e\}$. By Lemma~\ref{le1}, $|E_1|,|E_2|\geq 3$. For $i=1,2$, let $F_i$ be the subgraph of the link
graph of $v_i$ formed by the edges in $E_i$. If $\tau(F_i)\geq 3$, then  any edge $e\in E_{3-i}$ does not cover some edge
$f\in F_i$ and thus is disjoint from $f+v_1\in H$, a contradiction. Thus $\tau(F_1)\leq 2$
and $\tau(F_2)\leq 2$.

\smallskip
{\bf Case 1:} $\tau(F_1)=1$. 
Suppose $x_1$ is a dominating vertex in $F_1$. Since $|F_1|=|E_1|\geq 3$,
$x_1$ is {\bf the} dominating vertex in $F_i$ and
we may assume that $x_1x_2,x_1x_3,x_1x_4\in F_1$. But to cover these
$3$ edges, each edge in $F_2$ must contain $x_1$. Thus
 $H\subseteq H_0(n)$, as claimed.

\smallskip
{\bf Case 2:} $\tau(F_1)=\tau(F_2)=2$. If say $F_1$ contains a triangle $T=y_1y_2y_3$, then
$F_2$ cannot contain an edge not in $T$ and thus $F_2=T$ and by symmetry $F_1=T$. Thus
$H$ is contained in $H_4$. 

So the remaining case is that each of $F_i$ contains a matching
$M_i=\{z_{1,i}z_{1,i}',z_{2,i}z_{2,i}'\}$. Since  each edge of $F_1$ intersects each
edge of $F_2$, we may assume $z_{1,2}=z_{1,1},z'_{1,2}=z_{2,1},z_{2,2}=z'_{1,1},z'_{2,2}=z'_{2,1}$.
The only other edges that may have $F_2$ are $f_1=z_{1,1}z'_{2,1}$ and $f_2=z'_{1,1}z_{2,1}$. 
Since $|F_2|\geq 3$, we may assume $f_1\in F_2$.
Then the only third edge that $F_1$ may  contain is also $f_1$.
 It follows that $H$ is contained in $H_5$. This proves the main part of the theorem.

\smallskip
To prove part (a), assume  $H$ is an  intersecting $n$-vertex
$3$-graph with $|H|\geq 11$.
 Since $|K_5^3|=10<|H|$, $n\geq 6$. 
By Proposition~\ref{folk}, $\tau(H)\leq 2$. So part (a) is implied by the main claim of the theorem.
Part (b) follows from the fact that each of $H_3,H_4,H_5$ has $n+4$ edges.
 \qed

\section{Proof of Theorem~\ref{th4'}}
Let $H$ be as in the statement. By Theorem~\ref{fp}, $\tau(H)\leq 2$.
So, suppose a set $\{v_1,v_2\}$ covers all edges of $H$.
Let $E_0=\{e\in H\,:\,
 \{v_1,v_2\}\subset e\}$, and for $i=1,2$, let $E_i=\{e\in H\,:\,
v_{3-i}\notin e\}$. 

For $E_1\cup E_2$, construct the family $B(H)=B_1\cup B_2\cup\ldots B_r $  as in the previous proofs. Recall that by the minimality of the sets in $B_i$,
\begin{equation}\label{sper}
\mbox{$X \not\subseteq Y$ for all distinct $X, Y\in B(H)$,}
\end{equation}
and since $H$ is intersecting,
\begin{equation}\label{inter}
\mbox{$B(H)$ is intersecting.}
\end{equation}
If $B_1\neq \emptyset$, say $\{v_0\}\in B_1$, then by~\eqref{sper} and~\eqref{inter}, and
$B(H)=\{\{v_0\}\}$. This means either $H\subseteq H(n,r)$ (when $v_0\in \{v_1,v_2\}$), or
$H\subseteq H_0^r(n)$  (when $v_0\notin \{v_1,v_2\}$), and the theorem holds. 
So, let $B_1=\emptyset$.
 
Let $H'$ be obtained from $H$ by deleting all edges not containing a member of  $B'=B_2\cup B_3$. Then $|H-H'|\leq Cn^{r-4}$.
Since $\{v_1,v_2\}$ dominates $H$,
\begin{equation}\label{v1v2}
\mbox{each $D\in B'$ must contain either $v_1$ or $v_2$.}
\end{equation}
 For $i=1,2$, let $B'_i$ 
be the set of the members
of $B'$ containing $v_i$.

Define the auxiliary $3$-graph $H''$ with vertex set $V(H)$ as follows. The edges of $H''$ are all members of $B_3$ and each
triple $f$ that contains a member of $B_2$ and is contained in an $e\in H'$.

By~\eqref{inter}, $H''$ is intersecting. By~\eqref{v1v2},  $\tau(H'')\leq 2$.
If $\tau(H'')=1$, then $H'$ is a star. Suppose $\tau(H'')= 2$. By Theorem~\ref{th4},
$H''$ is contained in one of $H(n), H_0(n),  \ldots,H_5(n)$. But then
$H'$ is contained in one of $ H^r_0(n),  \ldots,H^r_5(n), EM(n,r,1)$, as claimed.
\qed

\section{Proof of Theorem~\ref{th1}}

 Recall that $r \ge 4, s\geq 1, n$ is sufficiently large and $H$ is an $n$-vertex  $r$-graph with $\nu(H)\leq s$ and $|H|>em(n,r,s-1)+hm''(n-s+1,r)$. 
We are to show that $V(H)$ contains a subset $Z=\{z_1,\ldots,z_{s-1}\}$ such that either $\tau(H-Z)=1$ or
$H-Z \subseteq HM(n-s+1,r,t)$ for some $t\in\{1,\ldots,r-1,n-s+1-r\}$ or $r=4$ and $H-Z \subseteq HM(n-s+1,4,0)$.

Define $B(H)$ and $B_i$ as in the previous proofs with the slight change that $T \in B(H)$ lies in an $(rs)^{|T|+1}$-sunflower (instead of an $(r+1)^{|T|}$-sunflower). 
Then the following claim holds (with an identical proof).

{\bf Claim.} $B_i$ contains no $(rs)^{i}$-sunflower.

Using the Claim and Lemma~\ref{er} we obtain $|B_i| < f((rs)^{i},i)$
for all $1 \le i \le r$.
As before,  setting $h=|B_1|$ we have 
$$ |H| \le \sum_{B \in B_1} {n-1 \choose r-1} + \sum_{i=2}^r \sum_{B \in B_i} {n-i \choose r-i}<h {n-1 \choose r-1} + (r-1)f((rs)^{r},r){n\choose r-2}.$$
Since $|H|>em(n,r,s-1)+ hm''(n-s+1, r) \sim s{n \choose r-1}$ and $n$ is large,  this immediately gives $h \ge s-1$. 
Consider distinct vertices $z_1, \ldots, z_{s-1} \in B_1$ and the set of edges  $F \subset H$ omitting $z_1, \ldots, z_{s-1}$.
 If $F$ is not intersecting, then let $e, e'$ be two disjoint edges in $F$. There exists a matching $e_1, \ldots, e_{s-1}$ in $H$ with 
$z_i \in e_i$ and $(e \cup e') \cap e_i = \emptyset$ for all $1 \le i \le s-1$. 
Note that we can produce the $e_i$ one by one since each $z_i$ forms the core of  an $(rs)^2$-sunflower
 in $H$ due to the definition of $B_1$.  We obtain the matching $e,e', e_1, \ldots, e_{s-1}$ contradicting $\nu(H) \le s$. 
Consequently, we may assume that $F$ is intersecting. Because
$|H|>em(n,r,s-1)+ hm''(n-s+1, r)$ we have
 $|F|>hm''(n-s+1, r)$.  Now we apply Theorem~\ref{th3} to $F$ to conclude that Theorem~\ref{th1} holds. \qed

\section{Concluding remarks}

Say that a hypergraph $H$ is $t$-{\em irreducible}, if
 $\nu(H)=t$ and $\nu(H-x)=t$ for every $x\in V(H)$.
Frankl~\cite{Frankl2013} presented a family of $n$-vertex $t$-irreducible $r$-graphs $PF(n,r,t)$ such that
 $$pf(n,r,t)=|PF(n,r,t)|\sim r{t-1\choose 2}{n\choose r-2}.$$
He also proved 

\begin{theorem}[\cite{Frankl2013}] \label{tpf}
Let $r \ge 4$, $t \ge 1$, and let $n$ be sufficiently large. Then  every 
 $n$-vertex $t$-irreducible $r$-graph $H$ has at most $pf(n,r,t)$ edges with equality
only if $H=PF(n,r,t)$.
\end{theorem}
Using this result, one can prove the following. 

\begin{lemma} \label{cor1int}
For every $r \ge 3$, $s\geq t\ge 2$, if $n$ is large, and $H$ is an $n$-vertex $r$-graph with  $\nu(H) = s$ and
$$|H|>em(n,r,s-t)+pf(n-s+t, r,t),$$ then there exists $X \subseteq V(H)$ with $|X|=s-t+1$ such that $\nu(H-X)=t-1$.
The bound on $|H|$ is sharp.
\end{lemma}
This in turn implies the following claim.

\begin{theorem} \label{cor2int}
For every $r \ge 3$ and $s\ge 2$ there exists $c>0$ such that the following holds. If $n$ is large, and $H$ is an $n$-vertex $r$-graph with  $\nu(H) = s$ and
$$|H|>em(n,r,s-2)+pf(n-s+2, r,2),$$ then either

1) there exists $H' \subset H$ with $|H'|<cn^{r-3}$ and $\tau(H-H')\le s$  or

2) there exist an $X \subset V(H)$ with $|X|=s-1$ and $u,v,w \in V(H-X)$ such that every edge of $H-X$  contains at least two elements of 
$\{u,v,w\}$.
\end{theorem}

We leave the details of the proofs to the reader.

Most of the proofs in this paper are rather simple applications of the early version of the Delta-system method. 
 There has been renewed interest in stability versions for problems in extremal set theory, 
so the general message of  this work is that the Delta-system method can quickly give some structural 
information about problems in extremal set theory, a fact that was already shown in several papers by Frankl and F\"uredi in the 1980's.
For more advanced recent applications of the Delta-system method, see the papers of F\"uredi~\cite{Furedi} and F\"uredi-Jiang~\cite{FJ}.
\\

{\bf Acknowledgment.} We thank Peter Frankl for helpful comments on an earlier version of the paper. 
We also thank Jozsef Balogh and Shagnik Das for attracting our attention to~\cite{HK}.

\end{document}